\documentclass[12pt]{amsart}
\usepackage{amsmath}
\usepackage{amsthm}
\usepackage{latexsym}
\usepackage{amssymb}
\DeclareMathOperator{\Char}{Char}

\DeclareMathOperator{\ad}{ad}

\DeclareMathOperator{\codim}{codim}
\newtheorem{theorem}{Theorem}
\newtheorem{proposition}{Proposition}
\newtheorem{lemma}{Lemma}
\newtheorem{corollary}{Corollary}

\newcounter{obsctr}

\renewcommand{\theequation}{\thesection.\arabic{equation}}
\begin{document}
\baselineskip 16pt
\def\A {{\mathcal{A}}}
\def\D {{\mathcal{D}}}
\def\R {{\mathbb{R}}}
\def\N {{\mathbb{N}}}
\def\C {{\mathbb{C}}}
\def\Z {{\mathbb{Z}}}
\def\l {\ell}
\def\ml {multline}
\def\multiline {\multline}
\def\lessim {\lesssim}
\def\d {\partial}
\def\phi{\varphi}
\def\epsilon{\varepsilon}
\title{Analytic Hypoellipticity at Non-Symplectic Poisson-Treves
  Strata for Sums of Squares of Vector Fields}    

\author{Antonio Bove}
\address{Dipartimento di Matematica, 
Universit\`a di Bologna, Piazza
di Porta San Donato 5, 40127 Bologna, Italy}
\email{Antonio.Bove@bo.infn.it} 
\author{David S. Tartakoff}
\address{Department of Mathematics, University
of Illinois at Chicago, m/c 249, 851 S.
Morgan St., Chicago IL  60607, USA}
\email{dst@uic.edu}
\date{\today}
\begin{abstract}
We consider an operator $ P $ which is a sum of squares of vector
fields with analytic coefficients. The operator has a non-symplectic
characteristic manifold, but the rank of the symplectic form $ \sigma
$ is not constant on $ \Char P $. Moreover the Hamilton foliation of
the non symplectic stratum of the Poisson-Treves stratification for $
P $ consists of closed curves in a ring-shaped open set around the
origin. We prove that then $ P $ is analytic hypoelliptic on that
open set. And we note explicitly that the local Gevrey hypoellipticity
for $ P $ is $ G^{k+1} $ and that this is sharp.
\end{abstract}
\maketitle
\pagestyle{myheadings}
\markboth{A. Bove and D. S. Tartakoff}{Analytic Hypoellipticity }
\section{Introduction}
\renewcommand{\theequation}{\thesection.\arabic{equation}}
\setcounter{equation}{0}
\setcounter{theorem}{0}
\setcounter{proposition}{0}  
\setcounter{lemma}{0}
\setcounter{corollary}{0} 
\setcounter{definition}{0}

The purpose of this paper is to study analytic hypoellipticity for
some sums of squares of vector fields having a non trivial
Poisson-Treves stratification. By this we mean that the stratification
has non trivial deep strata.

After Treves introduced the Poisson stratification associated to a sum
of squares of real analytic vector fields (\cite{Treves},
\cite{BoveTreves}), Hanges, \cite{Hanges}, observed that if we look at
analytic hypoellipticity in the sense of germs, there are operators
with a non symplectic characteristic set which are analytic
hypoelliptic in the sense of germs. Subsequently in \cite{BDT} a
class has been defined basically having the same properties of the
Hanges operator and being analytic hypoelliptic in the sense of
germs. It was also remarked that this fact is in no contradiction with
Treves' conjecture roughly stating that a sum of squares is analytic
hypoelliptic if and only if every stratum in its stratification is
symplectic. 

In all known examples we have analytic hypoellipticity in the sense
of germs when the characteristic manifold is non symplectic and is
actually a stratum of the stratification, which implies that the
symplectic form $ \sigma = d\xi \wedge dx $ has constant
rank. Moreover the canonical 1-form $ \omega = \xi dx $ does not
vanish. This implies that the characteristic manifold has particularly
simple microlocal models. 

Moreover the bicharacteristic curves, i.e. the Hamilton leaves of the
foliation, are closed curves foliating a given neighborhood of a
characteristic point on which hypoellipticity in the sense of germs
is obtained. 

In this paper we study two different cases. The first is that of a sum
of two squares for which the Poisson-Treves stratification has a
simplectic ``surface'' stratum and a deeper non-symplectic
stratum. The leaves of the Hamilton foliation are closed and foliate a
certain open subset of the stratum:
$$ 
P(t, x, D_{t}, D_{x})=
D_{t}^{2} + [x_{1} D_{2} - x_{2} D_{1} + t^{k} (x_{1} D_{1} + x_{2}
D_{2})]^{2},
$$
for $ k \geq 2 $ (see Section 2 for more details on its
stratification). 

It is well known that the above operator is $ G^{k+1} $ hypoelliptic
and not better when $ k \geq 2, $ and is analytic
hypoelliptic if $ k = 1 $.

We prove the following
\begin{theorem}
\label{th:1}
Let $ k \geq 2 $ and $ P $ be as above. Let $ U $ be an open subset of
$ \R^{3} $ in the variables $ (t, x_{1}, x_{2}) $ projecting on an
annulus of the form $ r_{1} < |x| < r_{2} $ and containing points
where $ t=0 $. Then $ P $ is analytic hypoelliptic (in the sense of
germs) at points in $ U \cap \Sigma_{2} $.
\end{theorem}

For the proof we need to use ideas, specifically localizations of high order derivatives adapted to the problem at hand and less straightforward than those of \cite{Tartakoff1}, \cite{Tartakoff2}, introduced by Derridj and Tartakoff in \cite{DT}.

Finally we discuss also another model operator which does not have
closed orbits; in this case the orbits foliate an annulus in the $ x
$-variables, but have $ \omega $- and $ \alpha $-limit sets that are
closed stationary orbits (see e.g. \cite{hartman}). 

In this case we have analytic regualrity of the solution if the closed
limit sets do not intersect the analytic wave front set of the
sulution. We do not attempt to prove any sort of analytic
hypoellipticity in the sense of germs in this case, since the fact
that the orbits are not closed does not seem to allow this.

\section{Proof of Theorem \ref{th:1}: Some preparations}
\renewcommand{\theequation}{\thesection.\arabic{equation}}
\setcounter{equation}{0}
\setcounter{theorem}{0}
\setcounter{proposition}{0}  
\setcounter{lemma}{0}
\setcounter{corollary}{0} 
\setcounter{definition}{0}

For $ k \in \N $, $ k \geq 2 $, let us consider the operator
\begin{equation}
\label{2.1}
P(t, x, D_{t}, D_{x})=
D_{t}^{2} + [x_{1} D_{2} - x_{2} D_{1} + t^{k} (x_{1} D_{1} + x_{2}
D_{2})]^{2}, 
\end{equation}
in the region $ 0 < r_{1} \leq r
\leq r_{2}, r=|x|, x= (x_{1}, x_{2}) \in \R^{2}.$

The Poisson-Treves stratification for $ P $ above is given by
\begin{eqnarray}
\label{2.2}
\Sigma_{1} & = & \{ \tau = 0, x_{1} \xi_{2} - x_{2}\xi_{1} + t^{k}
(x_{1} \xi_{1} + x_{2} \xi_{2})) = 0, t \neq 0 \}; \\
\Sigma_{2} & = & \{ \tau = 0, t= 0, x_{1}\xi_{2} - x_{2}\xi_{1} = 0,
\xi = (\xi_{1}, \xi_{2}) \neq 0 \}  \nonumber \\
\Sigma_{j} & = & \Sigma_{2}, \quad \text{for }  j \leq k,
\nonumber \\
\Sigma_{k+1} & = & \{0\}, \nonumber
\end{eqnarray}
where the last equation means that $ \Sigma_{k+1} $ is just the zero
section of $ R^{*}\R^{3} $.

We explicitly remark that $ \Sigma_{1} $ is a symplectic submanifold
of codimension 2, while $ \Sigma_{2} $ is not symplectic. Moreover $
\Sigma_{1} \cup \Sigma_{2} = \Char P $. 

Let us denote by
\begin{eqnarray}
\label{2.3}
X_{1} & = & D_{t}  \\
X_{2} & = & x_{1} D_{2} - x_{2} D_{1} + t^{k} (x_{1} D_{1} + x_{2} D_{2})\\
& = & D_\theta + t^kD_r = D_\theta + R
\nonumber 
\end{eqnarray}
so that
\begin{equation}
\label{2.4}
P = X_{1}^{2} + X_{2}^{2}.
\end{equation}
with 
\begin{equation} 
\label{basic brackets}
[X_1, X_2] = kt^{k-1}R, \text{ and } [R, X_j] = 0, j=1,2.
\end{equation}
We have the a priori estimate
\begin{equation}
\label{2.5}
\|v \|^{2}_{1/(k+1)} + \|X_{1} v\|^{2} + \|X_{2} v \|^{2} \leq C \left \{
|\langle Pv, v \rangle| + \|v \|^{2} \right\}, \qquad v\in C_0^\infty.
\end{equation}
where $ \| \cdot \| $ denotes the $ L^{2} $-norm in $
\R_{t} \times \R^{2}_{x} $.

\subsection{The general scheme}

Since the operator is subelliptic, the solution will be in $C^\infty .$ Additionally, since for $t\neq 0,$ the characteristic manifold of $P$ is symplectic, we know the solution is analytic for $t\neq 0.$ With a localizing function $\varphi (r)$ to be made precise below (but of Ehrenpreis type), and exploiting the maximality of the {\it a priori} estimate satisfied by $P,$ we will study $\|\varphi X^{p+1} u\|^2,$ each occurance of $X$ being $X_1$ or $X_2.$ Using the {\it a priori} estimate effectively will require moving one $X$ to the left of $\phi$ but this will not present a problem in the ensuing recursion.

We will immediately be led to estimate the bracket $|\langle[P,\varphi X^p]u, \varphi X^p u\rangle|.$
Upon iteration, using \ref{basic brackets}, we arrive, after at most $p$ iterations of the {\it a priori} estimate, to terms of the form 
\begin{equation}
\label{down by half}
\underline{C^p}\|(X)\varphi^{(p+1)}u\|^2 \text{ or } \underline{C^p}Cp!!\|(X)\varphi R^{p/2}u\|^2 
\end{equation}
and of course all the intermediate terms with some derivatives on $\varphi,$ some powers of $p,$ and some powers of $R,$ all with the generic bounds $$CC^pp^a\|(X)\varphi^{(b)}R^au\| \text{ with } p\sim b+2a.$$
Here $p!! = p(p-2)(p-4)\ldots,$ the value of $C$ may change from line to line but always independently of $u$ and the order of differentiation, and underlining a coefficient indicates the number of terms of the form which follows that occur. Finally, writing $(X)$ means that an $X=X_1$ or $X_2$ may or may not be present. 

When all $X$'s have been consumed in this way, we may no longer iterate effectively, and we must turn our attention to pure powers of $R,$ suitably localized. This will require a new localizing function and a construction we denote $(R^q)_\psi$ reminiscent of \cite{Tartakoff1} and \cite{Tartakoff2}, or more precisely \cite{DT1988}, which requires a special vector field $M$ which commutes especially well with both $X_1$ and $X_2,$ namely reproducing $X_1$ or generating $R.$ 

\subsection{The vector field $M$ and the localization}

We are fortunate to have a `good' vector field $M$ at our disposal which reproduces $R$ by bracketing with $X_2$: 
with 
\begin{equation}
\label{M}
M=\frac{t}{k}D_t,
\end{equation}
we have
\begin{equation}
\label{[M,X_2]}
[M, X_2] = R
\end{equation}

As localizing functions we shall use a nested family of Ehrenpreis-type functions as used by the second author in \cite{Tartakoff1}, \cite{Tartakoff2}. Given $N\in \N,$ the band $r\in (r_1, r_2)$ will contain $\log_2N$ nested subbands, $\Omega_k = \{r: r_{1_k}\leq r \leq r_{2_k}\}, r_{1_0}=r_1, r_{2_0}=r_2, k\leq \log_2N,$ with 
\begin{equation}
\label{d_k}
d_k=r_{1_k}-r_{1_{k-1}}=r_{2_k}-r_{2_{k-1}} = (r_2-r_1)\frac{1}{4k^2}
\end{equation} 
(so that $\sum d_k \leq r_2-r_1$)
and functions $\phi_k \equiv 1$ on $\Omega_{k}$ and supported in $\Omega_{k+1},$ such that with a constant $C=C_{r_2-r_1},$
\begin{equation}
\label{phi_k} |\phi_k^{(\ell)}(r)| \leq (C/d_k)^{\ell + 1} N_k^{\ell} \text{ for } \ell\leq N_k=N/{2^{k-1}}. 
\end{equation}
 The functions $\phi_k,$ but not the constant $C,$ depend on the choice of $N.$ In fact, we shall double the number of these functions, for technical reasons, $\phi_1, \tilde{\phi_1}, \phi_2, \tilde{\phi_2}, \ldots$ with $\phi_j$ and $\tilde{\phi_j}$ satisfying the same growth estimates. 
 
We note in passing, and will use later, that the growth estimates (\ref{phi_k}) imply the (weaker) growth estimates 
\begin{equation}
\label{phi_k'} |\phi_k^{(\ell)}(r)| \leq (C'/d_k)^{N_k^{\ell}+1} {\ell !} \text{ for } \ell\leq N_k=N/{2^{k-1}}. 
\end{equation}
(using the fact that $y^{1/y} \leq e$ for $y\geq 1.)$

Given the definition of $ M $ above, and for $ p \in \R $ large and $ j \in \{0, 1, \ldots, p\} $ we define the
expressions 
\begin{equation}
\label{N_j}
N_{j} = \sum_{j' = 0}^{j} a^{j}_{j'} \frac{M^{j'}}{j'!},
\end{equation}
where the $ a^{j}_{j'} $ denote rational numbers satisfying properties
that shall be made precise below which optimize commutation relations.

Finally we define our localizing operator, which is equal to $R^p$ where $\phi \equiv 1.$ We let
\begin{equation}
\label{2.9}
R^{p}_{\phi} = \sum_{j=0}^{p} \phi^{(j)} N_{j} R^{p-j}=\sum_{j=0}^{p} (R^j\phi) N_{j} R^{p-j}.
\end{equation}

\subsection{The Commutation Relations for $ R^{p}_{\phi} $}

For two vector fields $ Z $ and $ \tilde{Z} $ we shall frequently use
the formula 
\begin{equation}
\label{2.10}
[Z^{j}, \tilde{Z}] = \sum_{k=1}^{j} \binom{j}{k} \ad_{Z}^{k}
(\tilde{Z}) Z^{j-k},
\end{equation}
where
$$ 
\ad_{Z} (\tilde{Z}) = [Z, \tilde{Z}], \quad \ad_{Z}^{2}(\tilde{Z})
= [Z, [Z, \tilde{Z}] ]
$$
et cetera.

\subsection{The bracket $[X_2, R^p_\phi]$}

We first compute the commutator of $ X_2 $ with $ N_{j} $. We
have
$$
[ X_{2}, N_{j} ] = \sum_{j'=1}^{j} a^{j}_{j'} [ X_{2},
\frac{M^{j'}}{j'!}]
= - R^{k} \sum_{j'=1}^{j} a^{j}_{j'} \sum_{\ell=1}^{j'} \frac{1}{\ell
!} \frac{M^{j' - \ell}}{(j' - \ell)!} R,
$$
since 
$$ 
[ X_{2}, M ] = - t^{k} R.
$$

We seek to find coefficients $ a^{j}_{j'} $ so that
$$ 
[X_{2}, N_{j}] = - t^{k} N_{j-1} R,
$$
which will ensure that the bracket $[X_2, R^p_\phi]$ is free of the (poorly controlled) vector field $R$ (see below). Using 
 (\ref{2.8}), the necessary condition is that the $ a^{j}_{j'} $ must satisfy
\begin{equation}
\label{2.11}
\sum_{j'=1}^{j} \sum_{\ell=1}^{j'} a^{j}_{j'} \frac{1}{\ell !}
\frac{M^{j'-\ell}}{(j'-\ell)!} = \sum_{j_{1}=1}^{j-1} a^{j-1}_{j_{1}}
\frac{M^{j_{1}}}{j_{1}!}, 
\end{equation}
or
\begin{equation}
\label{2.12}
\sum_{s=1}^{j-\ell} a^{j}_{\ell+s} \frac{1}{s!} = a^{j-1}_{\ell},
\end{equation}
for $ \ell = 0, 1, \ldots, j-1 $.

We shall come back to condition (\ref{2.12}) later; for the time being
we may conclude the following
\begin{lemma}
\label{lemma:2.1}
With the coefficients $a^j_{j'}$ chosen as above, 
$$
[X_{2}, N_{j}] = - t^{k} N_{j-1} R,
$$
for every $ j \in \N $.
\end{lemma}
\begin{proposition}
\label{[X_2,Rpphi]}
$$[X_{2}, R^{p}_{\phi} ] = t^{k} \phi^{(p+1)} N_{p}.$$
\end{proposition}
\begin{proof}
Using the above Lemma we have, for these $a^j_{j'},$ 
$$[X_{2}, R^{p}_{\phi} ] = t^{k} \sum_{j=0}^{p} \phi^{(j+1)} N_{j}
R^{p-j} - t^{k} \sum_{j=1}^{p} \phi^{(j)} N_{j-1}R^{p+1-j}$$
$$= t^{k} \phi^{(p+1)} N_{p} + t^{k} \sum_{j=1}^{p} \phi^{(j)} N_{j-1}
R^{p-(j+1)} - t^{k} \sum_{j=1}^{p} \phi^{(j)} N_{j-1} R^{p+1-j}$$
$$= t^{k} \phi^{(p+1)} N_{p}.$$
\end{proof}

\subsection{The bracket $ [X_{1}, R^{p}_{\phi}]$}
First remark that
\begin{equation}
\label{2.14}
\ad_{M}^{\ell} (X_{1}) = \left( - \frac{1}{k}\right)^{\ell} X_{1},
\end{equation}
for every $ \ell \in N $. Therefore, 
$$ 
[X_{1}, N_{j}] = \sum_{j'=1}^{j} a^{j}_{j'} [X_{1},
\frac{M^{j'}}{j'!}] = - X_{1} \sum_{j'=1}^{j} \sum_{\ell=1}^{j'}
a^{j}_{j'} \left( - \frac{1}{k}\right)^{\ell} \frac{1}{\ell!}
\frac{M^{j'-\ell}}{(j'-\ell)!},
$$
so that we have
\begin{equation}
\label{2.15}
[X_{1}, R^{p}_{\phi}] = \sum_{j=1}^{p} \phi^{(j)} (- X_{1})
\sum_{j'=1}^{j} \sum_{\ell=1}^{j'} a^{j}_{j'}
\left(-\frac{1}{k}\right)^{\ell} \frac{1}{\ell!}
\frac{M^{j'-\ell}}{(j'-\ell)!} R^{p-j}.
\end{equation}

Our next goal is to prove the following lemma:
\begin{lemma}
\label{lemma:2.2}
For every $ j \in \N $ and $ \ell \in \{ 1, \ldots , j\} $ there exist
real constants $ \delta_{s} $, $ s=0, \ldots , j-2 $, such that
\begin{equation}
\label{2.16}
\sum_{h=1}^{j-\ell} a^{j}_{\ell+h} \left(-\frac{1}{k}\right)^{h}
\frac{1}{h!} = \sum_{h=1}^{j-\ell} \delta_{j-\ell -h} a^{\ell+h-1}_{\ell}.
\end{equation}
\end{lemma}
The above Lemma has an easy consequence:
\begin{lemma}
\label{lemma:2.3}
For every $ j \in \N $  there exist
real constants $ \gamma_{s} $, $ s=0, \ldots , j-1 $, such that
\begin{equation}
\label{2.17}
\sum_{j'=1}^{j} \sum_{\ell=1}^{j'}
a^{j}_{j'}\left(-\frac{1}{k}\right)^{\ell}  \frac{1}{\ell!}
\frac{M^{j'-\ell}}{(j'-\ell)!} = \sum_{s=0}^{j-1} \gamma_{j-s} N_{s}.
\end{equation}
\end{lemma}
\begin{proof}[Proof of Lemma \ref{lemma:2.3}]
The identity (\ref{2.17}) can be restated as
$$ 
\sum_{s=0}^{j-1} \sum_{j'=s+1}^{j} a^{j}_{j'}
\left(-\frac{1}{k}\right)^{j'-s} \frac{1}{(j'-s)!} \frac{M^{s}}{s!} =
\sum_{s=0}^{j-1} \sum_{h=0}^{s} \gamma_{j-s} a^{s}_{h} \frac{M^{h}}{h!}
$$
or
$$ 
\sum_{s=0}^{j-1} \sum_{j'=s+1}^{j} a^{j}_{j'}
\left(-\frac{1}{k}\right)^{j'-s} \frac{1}{(j'-s)!} \frac{M^{s}}{s!} =
\sum_{\ell=0}^{j-1}\sum_{s=\ell}^{j-1} \gamma_{j-s} a^{s}_{\ell}
\frac{M^{\ell}}{\ell!}. 
$$
From which we get
$$ 
\sum_{j'=\ell+1}^{j} a^{j}_{j'} \left(-\frac{1}{k}\right)^{j'-\ell}
\frac{1}{(j'-\ell)!} = \sum_{s=\ell}^{j-1} \gamma_{j-s} a^{s}_{\ell}.
$$
Now the latter identity can be rewritten as
$$ 
\sum_{h=1}^{j-\ell} a^{j}_{\ell+h} \left(-\frac{1}{k}\right)^{h}
\frac{1}{h!} = \sum_{h=1}^{j-\ell} \delta_{j-\ell -h} a^{\ell+h-1}_{\ell},
$$
for any $ \ell = 1, \ldots, j-1 $, and this is in the statement of
Lemma \ref{lemma:2.2}.
\end{proof}
\begin{proof}[Proof of Lemma \ref{lemma:2.2}]
In order to prove Lemma \ref{lemma:2.2} we must analyse the recurrence
relation (\ref{2.12}):
$$
\sum_{h=1}^{j-\ell} a^{j}_{\ell+h} \frac{1}{h!} = a^{j-1}_{\ell},
$$
for $ \ell = 0, 1, \ldots, j-1 $.

Another way of rewriting the above relation is the following:
\begin{equation}
\label{2.18}
\begin{bmatrix}
 1 & \frac{1}{2!} & \cdots & \frac{1}{(j-1)!} & \frac{1}{j!} \\[7pt] 
 0 & 1 & \cdots & \frac{1}{(j-2)!} & \frac{1}{(j-1)!} \\[7pt] 
 \vdots & \vdots & \ddots & \vdots & \vdots \\[7pt] 
 0 & 0 & \cdots & 1 & \frac{1}{2!} \\[7pt] 
 0 & 0 & \cdots & 0 & 1 \\
\end{bmatrix}
\begin{bmatrix}
  a^{j}_{1} \\[7pt]
   \vdots \\[7pt]
  a^{j}_{j}
\end{bmatrix}
=
\begin{bmatrix}
  a^{j-1}_{0} \\[7pt]
  \vdots \\[7pt]
  a^{j-1}_{j-1}
\end{bmatrix}
.
\end{equation}
Note that on the left hand side there are no terms of the form $
a^{j}_{0}
$, which means that we are free to choose those coefficients. We shall
choose $ a^{0}_{0} = 1 $ for the sake of simplicity, leaving the
others undetermined.

We point out that the matrix in the above formula is clearly
invertible and that it can be written as
$$ 
I_{j} + \frac{1}{2!} J_{j} + \frac{1}{3!} J_{j}^{2} + \cdots +
\frac{1}{j!} J_{j}^{j-1} = \int_{0}^{1} e^{t J_{j}} dt,
$$
where $ J_{j} $ denotes the standard $ j \times j $ Jordan matrix
$$ 
J_{j} =
\begin{bmatrix}
 0 & 1 & 0 & \cdots & 0 & 0 \\
 0 & 0 & 1 & \cdots & 0 & 0 \\ 
 \vdots & \vdots & \vdots & \ddots & \vdots & \vdots \\
 0 & 0 & 0 & \cdots & 1 & 0 \\
 0 & 0 & 0 & \cdots & 0 & 1 \\ 
 0 & 0 & 0 & \cdots & 0 & 0 \\ 
\end{bmatrix}.
$$
Using, for example, formula (\ref{2.18}) we may easily see that,
inverting the matrix, we obtain
\begin{equation}
\label{2.19}
\begin{bmatrix}
  a^{j}_{1} \\[7pt]
   \vdots \\[7pt]
  a^{j}_{j}
\end{bmatrix} =
\begin{bmatrix}
 c_{0} & c_{1} & \cdots & c_{j-2} & c_{j-1} \\[7pt] 
 0 & c_{0} & \cdots & c_{j-3} & c_{j-2} \\[7pt] 
 \vdots & \vdots & \ddots & \vdots & \vdots \\[7pt] 
 0 & 0 & \cdots & c_{0} & c_{1} \\[7pt] 
 0 & 0 & \cdots & 0 & c_{0} \\
\end{bmatrix}
\begin{bmatrix}
  a^{j-1}_{0} \\[7pt]
  \vdots \\[7pt]
  a^{j-1}_{j-1}
\end{bmatrix},
\end{equation}
where $ c_{0}=1 $, $ c_{1} = - \frac{1}{2!} $ and the other $ c_{m} $
can be computed by a triangular relation. In particular, using the
structure of the matrix, we obtain that
\begin{equation}
\label{2.20}
\begin{bmatrix}
  a^{j}_{\ell+1} \\[7pt]
   \vdots \\[7pt]
  a^{j}_{j}
\end{bmatrix} =
\begin{bmatrix}
 c_{0} & c_{1} & \cdots & c_{j-\ell-2} & c_{j-\ell-1} \\[7pt] 
 0 & c_{0} & \cdots & c_{j-\ell-3} & c_{j-\ell-2} \\[7pt] 
 \vdots & \vdots & \ddots & \vdots & \vdots \\[7pt] 
 0 & 0 & \cdots & c_{0} & c_{1} \\[7pt] 
 0 & 0 & \cdots & 0 & c_{0} \\
\end{bmatrix}
\begin{bmatrix}
  a^{j-1}_{\ell} \\[7pt]
  \vdots \\[7pt]
  a^{j-1}_{j-1}
\end{bmatrix},
\end{equation}
for $ \ell = 0, 1, \ldots, j-1 $.

Another way of writing the above identity is
\begin{equation}
\label{2.21}
a^{j}_{\ell+h} = \sum_{s=h}^{j-\ell} c_{s-h} a^{j-1}_{\ell -1 +s},
\end{equation}
for $ h = 1, 2, \ldots, j-\ell $.

Iterating, we get
$$ 
a^{j-t}_{\ell-t+h} = \sum_{s=h}^{j-\ell} c_{s-h} a^{j-t-1}_{\ell-t-1+s},
$$
for $ t = 0, 1, \ldots , j-\ell-1 $.

Let us now fix an $ h \in \{1, \ldots, j-\ell\} $. Then we have
\begin{eqnarray*}
a^{j}_{\ell+h} & = & \sum_{s_{1}=h}^{j-\ell} c_{s_{1}-h}
a^{j-1}_{\ell-1+s} \\
   & = & \sum_{s_{1}=h}^{j-\ell} \sum_{s_{2}=s_{1}}^{j-\ell}
   c_{s_{1}-h} c_{s_{2}-s_{1}} a^{j-2}_{\ell-2+s_{2}} \\
   & = & \sum_{s_{1}=h}^{j-\ell} \sum_{s_{2}=s_{1}}^{j-\ell} \cdots
   \sum_{s_{h} = s_{h-1}}^{j-\ell} c_{s_{1}-h} c_{s_{2}-s_{1}} \cdots
   c_{s_{h} - s_{h-1}} a^{j-h}_{\ell-h+s_{h}}
\end{eqnarray*}
The latter sum can be written as
\begin{eqnarray*}
a^{j}_{\ell+h} & = & c_{0}^{h} a^{j-h}_{\ell} +
\sum_{s_{1}=h}^{j-\ell} \sum_{s_{2}=s_{1}}^{j-\ell} \cdots\sum_{\substack{s_{h}
= s_{h-1}\\s_{h} > h}}^{j-\ell} c_{s_{1}-h} c_{s_{2}-s_{1}} \cdots
   c_{s_{h} - s_{h-1}} a^{j-h}_{\ell-h+s_{h}}\\
& = &  c_{0}^{h} a^{j-h}_{\ell} +
\sum_{s_{1}=h}^{j-\ell} \sum_{s_{2}=s_{1}}^{j-\ell} \cdots\sum_{\substack{s_{h}
= s_{h-1}\\s_{h} > h}}^{j-\ell} \sum_{s_{h+1}=s_{h}}^{j-\ell}
c_{s_{1}-h} \cdots c_{s_{h+1}-s_{h}} a^{j-h-1}_{\ell-h-1+s_{h+1}}
\end{eqnarray*}
The latter sum allows us to compute the coefficient of $
a^{j-h-1}_{\ell} $, by picking all terms for which one of the $
s_{j} $ is equal to $ h+1 $, for $ j = 1, 2, \ldots , h+1 $.

Iterating this procedure, i.e. using the recursion relation until we
obtain a coefficient $ a^{*}_{*} $ where the lower index is equal to $
\ell $, we may express the coefficient $ a^{j}_{\ell+h} $ as a linear
combination of $ a^{*}_{\ell} $; the above formulas show that we may
actually write
\begin{equation}
\label{2.22}
a^{j}_{\ell+h} = \sum_{\sigma=0}^{j-\ell-h} \alpha_{j-\ell-\sigma}
a^{\ell+\sigma}_{\ell}, 
\end{equation}
for $ h = 1, 2, \ldots , j-\ell $.

We point out explicitly that up to this point we have only used the
recurrence relation (\ref{2.12}). Let us now denote by $ A_{h} $ the
collection of real numbers $ A_{h} = (- k^{-1})^{h} h!^{-1}
$. Then it is evident that
$$ 
\sum_{h=1}^{j-\ell} a^{j}_{\ell+h} A_{h} = \sum_{\sigma=0}^{j-\ell-1}
\delta_{j-\ell-\sigma} a^{\ell+\sigma}_{\ell},
$$
where $ \delta_{j-\ell-\sigma} = \sum_{h=1}^{j-\ell-\sigma} A_{h} $,
and this is the statement of the Lemma.
\end{proof}
\begin{lemma}[See \cite{DT} and \cite{hirzebruch}]
\label{lemma:2.4}
Let us consider the recurrence relation (\ref{2.12}):
$$
\sum_{s=1}^{j-\ell} a^{j}_{\ell+s} \frac{1}{s!} = a^{j-1}_{\ell}.
$$
Setting
\begin{equation}
\label{aa}
a_{\ell}^{j} = \frac{1}{(j-\ell)!}  \left( \left[ \frac{t}{e^{t} - 1}
  \right]^{j+1}\right)^{(j-\ell)} (0),
\end{equation}
we obtain a solution of the above recurrence satisfying the boundary
conditions $ a_{j}^{j} = 1 $ and $ a^{j}_{0} = (-1)^{j} $, $ j \geq 0
$. Moreover this is the only power series with rational coefficients
satisfying (\ref{2.12}) and the above boundary conditions.
\end{lemma}
\begin{proof}
By a simple computation we have
\begin{eqnarray*}
a_{\ell}^{j-1} & = & \frac{1}{(j-1-\ell)!} \left(
  \left(\frac{t}{e^{t}-1}\right)^{j}\right)^{(j-1-\ell)}(0) \\
& = & \frac{1}{(j-1-\ell)!} \left(
  \left(\frac{t}{e^{t}-1}\right)^{j+1} \frac{e^{t}-1}{t}
\right)^{(j-1-\ell)}(0) \\
& = & \sum_{h=0}^{j-1-\ell} \frac{1}{(j-1-\ell- h)!} \left(
  \left(\frac{t}{e^{t}-1}\right)^{j+1}\right)^{(j-1-\ell-h)}(0)
\frac{1}{(h+1)!}  \\
& = & \sum_{p=1}^{j-\ell} \frac{1}{(j-1-p)!} \left(
  \left(\frac{t}{e^{t}-1}\right)^{j+1}\right)^{(j-\ell-p)}(0)
\frac{1}{p!} \\
& = & \sum_{p=1}^{j-\ell} a_{\ell+p}^{j} \frac{1}{p!}.
\end{eqnarray*}
Where we used the fact that 
$$ 
\frac{1}{h!} \left( \frac{e^{t}-1}{t} \right)^{(h)}(0) = \frac{1}{(h+1)!}.
$$
Moreover we have $ a_{j}^{j} = 1 $ for every $ j \geq 0  $. As for the
other boundary condition, first we remark that
$$ 
a^{j}_{0} = \frac{1}{j!}
\left(\left[\frac{t}{e^{t}-1}\right]^{j+1}\right)^{(j)} (0),
$$
i.e. $ a^{j}_{0} $ is the coefficient of $ t^{j} $ in the power series
of $ Q(t) $,
$$ 
Q(t) = \left( \frac{t}{e^{t}-1} \right)^{j+1}.
$$
Thus 
$$ 
a^{j}_{0} = \frac{1}{2 i \pi} \int_{\gamma}
\left(\frac{1}{e^{z}-1}\right)^{j+1} dz,
$$
where $ \gamma $ is a smooth curve encircling the origin in $ \C $.

Changing variables $ w = e^{z} -1 $, so that the origin is mapped to
the origin and $ \gamma $ is mapped to another smooth curve encircling
the origin that we still denote by $ \gamma $, we have
$$ 
a^{j}_{0} = \frac{1}{2 i \pi} \int_{\gamma} w^{-(j+1)} (w + 1)^{-1} dz
= (-1)^{j}.
$$
The uniqueness is proved in \cite{hirzebruch}. This end the proof of
the lemma.
\end{proof}

As a consequence of the preceding Lemmas we may now state the
\begin{proposition}
\label{[X_1,Rpphi]}
The commutator of $ X_{1} $ with the localizing operator $
R^{p}_{\phi} $ has the form
\begin{equation}
\label{2.23}
[X_{1}, R^{p}_{\phi}] = -X_{1} \sum_{\ell=0}^{p-1} \delta_{\ell}
R^{p-\ell-1}_{\phi^{(\ell+1)}}, 
\end{equation}
where $ \phi^{(j)} $ denotes the $ j $-th derivative $
(r\partial_{r})^{j} \phi$ and $\delta_\ell = \sum_1^\ell\frac{1}{k^hh!} \leq 1.$
\end{proposition}

From Lemma \ref{lemma:2.4} we have the
\begin{corollary}
\label{cor:2.1}
For every $ j \geq 0 $ and $ \ell \in \{0, \ldots, j\} $ we have
$$ 
| a^{j}_{\ell} | \leq c^{j},
$$
for a suitable universal positive constant $ c $.
\end{corollary}
\begin{proof}
From (\ref{aa}) we have that
$$ 
a^{j}_{\ell} = \frac{1}{2 i \pi} \int_\gamma
\left(\frac{t}{e^{t}-1}\right)^{j+1} t^{-(j-\ell+1)} \ dt,
$$
where $ \gamma $ is a circle of fixed radius around the origin. Since
the function under the integral sign may be estimated by a positive
constant (depending on the radius of $ \gamma $) raised to the power $
j $, the corollary follows.
\end{proof}

\section{Proof of Theorem \ref{th:1}}
\renewcommand{\theequation}{\thesection.\arabic{equation}}
\setcounter{equation}{0}
\setcounter{theorem}{0}
\setcounter{proposition}{0}  
\setcounter{lemma}{0}
\setcounter{corollary}{0} 
\setcounter{definition}{0}

In this section we prove that $ P $ is analytic hypoelliptic
in any open set of the form $ \Omega = \{(t, x) \in \R^{3}\ |\ r_{1} <
|x| < r_{2}, t \in (-\delta, \delta) \}, \delta > 0.$

The maximal estimate may be restated to allow $X$ to appear to the right or left of the localizing function (where $\psi' = X\psi$): 

\begin{equation}
\label{apriori balanced}
\|\psi X^p u\|^2_{\frac{1}{k}} + \|X\psi X^pu\|^2+\|\psi X^{p+1}u\|^2 \lesssim 
\end{equation}
$$\lesssim |\langle P\psi X^pu, \psi X^pu\rangle|+\|\psi'X^p u\|^2 \lesssim
$$
$$\lesssim |\langle \psi X^p Pu, \psi X^pu\rangle|+ |\langle [X^2,\psi X^p]u, \psi X^p u\rangle|+\|\psi'X^p u\|^2 $$

Now
$$|\langle[X^2,\psi X^p]u, \psi X^p u\rangle| \leq  |\langle[X,\psi X^p]u, X\psi X^p u\rangle|
+|\langle[X,\psi X^p]Xu, \psi X^p u\rangle|$$
$$ \lesssim |\langle \psi' X^p u, X\psi X^p u\rangle| +  \underline{p}|\langle t^{k-1}\psi RX^{p-1}u, X\psi X^p u\rangle|$$
$$+ |\langle\psi X^pXu, \psi' X^pu\rangle| + \underline{p}|\langle t^{k-1}\psi X^pRu, \psi X^pu\rangle|
$$
$$\leq \epsilon \left\{ \|X\psi X^pu\|^2+\|\psi X^{p+1}u\|^2\right\} + C_\epsilon \left\{\| \psi' X^p u\|^2 +  (\underline{p}\|\psi X^{p-1}Ru\|)^2\right\}$$
where we have freely exchanged $\psi$ and $\psi'$ on the two sides of the inner product when no derivatives intervened. Note that $[X,R]=0.$

In all, 
$$\|\psi X^p u\|^2_{\frac{1}{k}} + \|X\psi X^pu\|^2+\|\psi X^{p+1}u\|^2 $$
$$\lesssim 
\|\psi X^p Pu\|^2 +\|\psi' X^p u\|^2 +  (\underline{p}\|\psi RX^{p-1}u\|)^2 $$

Iterating this inequality until there remain no $X$'s on the right, 

\begin{equation}
\label{XtoR}\|\psi X^p u\|^2_{\frac{1}{k}} + \|X\psi X^pu\|_{L^2}^2+\|\psi X^{p+1}u\|_{L^2}^2 
\end{equation}
$$\leq C^p\left\{\sup_{j+2d\leq p}\{p^d\|\psi^{(j)} R^dX^{p-j-2d} Pu\|_{L^2}^2\} +\sup_{j+2d= p}
(p^d\|\psi^{(j)} R^d u\|_{L^2})^2 \right\}. $$

The first term on the right can be estimated directly (even taken to be zero, using the Cauchy-Kowalevska theorem). For the second, we will take the localizing function out of the norm and introduce one of the $(R^d)_{\tilde{\psi}} \equiv R^d$ on the support of $\psi.$ Thus for such $\tilde{\psi},$ and taking $Pu=0$ for simplicity,

\begin{equation}
\|\psi X^p u\|_{\frac{1}{k}} + \|X\psi X^pu\|_{L^2}+\|\psi X^{p+1}u\|_{L^2} \leq 
\sup_{j+2d= p}
p^d\sup|\psi^{(j)}|\| (R^d)_{\tilde{\psi}} u\|_{L^2}
\end{equation}

 For convenience we recall the bracket relations and the few important definitions (for generic $\phi$):
$$[X_{1}, R^{b}_{\phi}] = -X_{1} \sum_{\ell=0}^{b-1} \delta_{\ell}
R^{b-\ell-1}_{\phi^{(\ell+1)}}, \quad |\delta_\ell |\leq 1$$
$$[X_{2}, R^{b}_{\phi} ] = t^{k} \phi^{(b+1)} N_{b}.$$
$$N_{b} = \sum_{b' = 0}^{b} a^{b}_{b'} \frac{M^{b'}}{b'!}, \qquad M=\frac{t}{k}D_t, \qquad |a^b_{b'}|\leq c^b.$$

As above, we use the {\it a priori} estimate, but now on $v=(R^d)_\phi u:$
\begin{equation}
\label{aprioriRdphi}
\|(R^d)_\phi u\|^2_{\frac{1}{k}} + \|X(R^d)_\phi u\|^2+\| (R^d)_\phi Xu\|^2 \lesssim 
\end{equation}
$$\lesssim |\langle P(R^d)_\phi u, \psi (R^d)_\phi u\rangle|+\|[X,(R^d)_\phi] u\|^2 \lesssim
$$
$$\lesssim |\langle (R^d)_\phi Pu, (R^d)_\phi u\rangle|+ |\langle [X^2,(R^d)_\phi]u, (R^d)_\phi u\rangle|+\|[X,(R^d)_\phi] u\|^2. $$
Again, taking $Pu=0,$ and expanding $[X^2,(R^d)_\phi]=X[X,(R^d)_\phi]+[X,(R^d)_\phi]X,$ we find, as before, with a weighted Schwarz inequality and integrating by parts one $X=-X^*,$ 
\begin{equation}
\label{RdphiX}\|(R^d)_\phi u\|^2_{\frac{1}{k}} + \|X(R^d)_\phi u\|^2+\| (R^d)_\phi Xu\|^2
\end{equation}
$$\lesssim |\langle [X,(R^d)_\phi]Xu, (R^d)_\phi u\rangle|+ \|[X,(R^d)_\phi] u\|^2
$$

Now on the right, when $X=X_1,$ the result, as we saw above, still has an $X_1,$ which we integrate by parts in the case of the inner product:
\begin{equation}\label{RdphiX1}|\langle [X_1,(R^d)_\phi]X_1u, (R^d)_\phi u\rangle|+ \|[X_1,(R^d)_\phi] u\|^2\end{equation}
$$\leq \epsilon\|X_1(R^d)_{\phi}u\|^2+C_\epsilon\sum_{d_1=1}^d \|(R^{d-d_1})_{\phi^{(d_1)}}X_1 u\|^2$$
 On the other hand, when $X=X_2,$ we have nearly pure powers of $tD_t,$ 
which it will be necessary to convert into pure powers of $X_1=D_t$ (from which we started, but, we note, of at most half the order). 

\begin{proposition} $$(tD_t)^j =  \sum_{\ell=1}^{j} B^j_\ell\;t^\ell D_t^\ell
$$
where
$$B^j_\ell =\sum_{m=0}^{\ell-1 }  \frac{(-1)^m(\ell -m)^{j-1}}{m!(\ell -m-1)!}=\sum_{m=0}^{\ell-1 }  \frac{(-1)^m(\ell -m)^{j}}{m!(\ell -m)!}
$$
so that for all $v,$ pointwise,
$$\frac{|(tD_t)^j v|}{j!} \leq C^j \sum_{\ell=1}^j \frac{|t^\ell D_t^\ell v|}{\ell !}$$
and hence in $|t|<1,$ 
$$|N_{b}v| = \left|\sum_{b' = 0}^{b} a^{b}_{b'} \frac{M^{b'}}{b'!}v\right|= \left|\sum_{b' = 0}^{b} (\frac{1}{k})^{b'}a^{b}_{b'} \frac{(tD_t)^{b'}}{b'!}v\right|\leq C^b\sup_{b'\leq b}\left|\frac{X_1^{b'}v}{b'!}\right|$$
\end{proposition}

The particular expression for the coefficients $B^j_\ell$ is proved by induction and can be understood by a kind of over-counting/under- counting argument. 

Thus for $X_2,$
\begin{equation}\label{RdphiX2}|\langle [X_2,(R^d)_\phi]X_2u, (R^d)_\phi u\rangle|+ \|[X_2,(R^d)_\phi] u\|^2\end{equation}
$$\leq |\langle t^k\phi^{(d+1)}N_d X_2u, (R^d)_\phi u\rangle|+ \|t^k\phi^{(d+1)}N_d u\|^2 $$
$$\leq C^d \sup_{d'\leq d}\|\phi^{(d+1)}\frac{X_1^{d'}(X_2)u}{d'!}\|^2 + \epsilon \|(R^d)_\phi u\|^2$$
or in all, 
\begin{equation}
\label{aprioriRdphiall}
\|(R^d)_\phi u\|^2_{\frac{1}{k}} + \|X(R^d)_\phi u\|^2+\| (R^d)_\phi Xu\|^2 
\end{equation}
$$\leq C\sum_{d_1=1}^d \|(R^{d-d_1})_{\phi^{(d_1)}}X_1 u\|^2+C^d \sup_{d'\leq d}\|\phi^{(d+1)}\frac{X_1^{d'}(X_2)u}{d'!}\|^2$$
(i.e., with or without $X_2$ in the last term). Iterating on the first term on the right, 
eventually only the last term survives:
\begin{equation}
\label{aprioriRdphiallb}
\|(R^d)_\phi u\|^2_{\frac{1}{k}} + \|X(R^d)_\phi u\|^2+\| (R^d)_\phi Xu\|^2\end{equation} 
$$\leq C^d \sup|\phi^{(d+1)}|\sup_{d'\leq d}\|\tilde{\phi}\frac{X^{d'+1}u}{(d'+1!)}\|^2$$
for any $\phi, \tilde{\phi}$ with $\tilde{\phi}\equiv 1$ on the support of $\phi.$ 

Recalling the previous bound 
\begin{equation}
\|\psi X^p u\|^2_{\frac{1}{k}} + \|X\psi X^pu\|^2_{L^2}+\|\psi X^{p+1}u\|^2_{L^2}\end{equation}
$$ \leq 
\sup_{j+2d= p}
p^d\sup|\psi^{(j)}|\| (R^d)_{\tilde{\psi}} u\|^2_{L^2},
$$
valid for any $\tilde{\psi}\equiv 1$ on the support of $\psi,$ we have the choice of starting with $X$'s, reducing the order by half, introducing $(R^d)_\phi$ and iterating that until we are back to $X$'s or start with $(R^d)_\phi,$ reduce to $X$'s until they bracket to yield pure $R$'s at half the order. In either order, after one full cycle, we need a new localizing function each time $(R^d)_\phi$ is put together. Thus in starting with $N$ derivatives to estimate, after $\log_2N$ full cycles, the number of free derivatives on $u$ will be only a bounded number. 

For definiteness, we follow the cycle starting with powers of $X$'s, and introduce for a moment the new norms
$$|||\psi, X^p, u||| = \|\psi X^p u\|_{\frac{1}{k}} + \|X\psi X^pu\|_{L^2}+\|\psi X^{p+1}u\|_{L^2}$$
and
$$|||(R^d)_\phi,u||| = \|(R^d)_\phi u\|_{\frac{1}{k}} + \|X(R^d)_\phi u\|_{L^2}+\| (R^d)_\phi Xu\|_{L^2},$$
so that the above may be written
\begin{equation}
|||(R^d)_\phi,u||| 
\leq C^d \sup_{d'\leq d}\frac{1}{(d'+1)!}|||\phi^{(d+1)},X^{d'+1},u|||\end{equation}
for any $\phi$ and 
\begin{equation}
|||\psi, X^p, u|||  \leq \sup_{j+2d= p} p^d\sup |\psi^{(j)}|\,||| (R^d)_{\tilde{\psi}},u|||,
\end{equation}

Thus we start with $\psi=\phi_1$ (the first in the sequence of precisely nested localizing functions (cf. (\ref{phi_k})) for a fixed $N=N_1\in\N$:
$$\frac{\|X^{N_1}u\|_{L^2(\Omega_1)}}{N_1!}\leq \frac{|||\phi_1, X^{N_1}, u|||}{N_1!}\leq \sup_{N_2\leq \widetilde{N_1}\leq N_1}\frac{|||\phi_1, X^{\widetilde{N_1}}, u|||}{\widetilde{N_1}!}$$
$$ \leq \sup_{{\ell_1+2\delta_1= {\widetilde{N_1}}} \atop N_2\leq \widetilde{N_1}\leq N_1} \frac{{\widetilde{N_1}}^{\delta_1}\sup |\phi_1^{(\ell_1)}|\,||| (R^{\delta_1})_{\tilde{\phi}},u|||}{\widetilde{N_1}!}$$
$$\leq \sup_{{\ell_1+2\delta_1= {\widetilde{N_1}}} \atop N_2\leq \widetilde{N_1}\leq N_1}  \frac{\left(\frac{C}{d_1}\right)^\ell \widetilde{N_1}^{\ell_1+\delta_1}||| (R^{\delta_1})_{\tilde{\phi}},u|||}{\widetilde{N_1}!}$$
with any $\tilde{\phi}\equiv 1 $ near the support of $\phi_1.$

Now there is some freedom in the choice of $\tilde{\phi},$ since all that we have required is that it be one on the support of $\phi_1,$ and we pick the largest index $k$ consistent with $\delta_1,$ i.e., $N_{k}\geq \delta_1 \geq N_{k+1}$ and $k\geq 2$ since $\ell_1+2\delta_1= N_1.$)
Thus with $\tilde{\phi}=\phi_k$ and 
together with the other estimate:
\begin{equation}
|||(R^\delta)_{\phi_k},u||| 
\leq C^\delta \sup |\phi_k^{(\delta+1)}|\sup_{\delta'\leq \delta}\frac{|||\widetilde{\phi_k},X^{\delta'+1},u|||}{(\delta'+1)!},\end{equation}
we arrive at 
$$\frac{\|X^{{N_1}}u\|^2_{L^2(\Omega_1)}}{{{N_1}}^{{N_1}}}\leq \sup_{N_2\leq \widetilde{N_1}\leq N_1}\frac{|||\widetilde{\phi_1}, X^{\widetilde{N_1}}, u|||}{{\widetilde{N_1}}^{\widetilde{N_1}}}$$
$$\leq \sup_{{\ell+2\delta= {\widetilde{N_1}}}\atop {N_2\leq \widetilde{N_1}\leq N_1+1}}  \frac{\left(\frac{C}{d_1}\right)^\ell \widetilde{N_1}^{\ell+\delta}C^\delta \sup|\phi_k^{(\delta+1)}|}{{\widetilde{N_1}}^{\widetilde{N_1}}}\sup_{{N_{k+1}\leq \widetilde{N_k}\leq N_k+1}\atop k\geq 2}\frac{|||\widetilde{\phi_k},X^{\widetilde{N_k}},u|||}{\widetilde{N_k}^{\widetilde{N_k}}}$$
$$\leq \sup_{{\ell+2\delta= {\widetilde{N_1}}}\atop {N_2\leq \widetilde{N_1}\leq N_1}} \frac{ \left(\frac{C}{d_1}\right)^\ell \widetilde{N_1}^{\ell+\delta}C^\delta \left(\frac{C}{d_k}\right)^{\delta+1}N_k^{(\delta+1)}}{{\widetilde{N_1}}^{\widetilde{N_1}}}\sup_{{N_{k+1}\leq \widetilde{N_k}\leq N_k+1}\atop k\geq 2}\frac{|||\widetilde{\phi_k},X^{\widetilde{N_k}},u|||}{\widetilde{N_k}^{\widetilde{N_k}}}$$
$$\leq C^{\widetilde{N_1}}\sup_{\ell+2\delta= {\widetilde{N_1}}}\frac{d_1^{-\ell}d_k^{-(\delta+1)}}{2^{\ell+\delta}2^{k(\delta+1)}}\sup_{{N_{k+1}\leq \widetilde{N_k}\leq N_k+1}\atop k\geq 2}\frac{|||\widetilde{\phi_k},X^{\widetilde{N_k}},u|||}{\widetilde{N_k}^{\widetilde{N_k}}}.$$

Now the expressions in the first supremum increase as $k$ decreases, bounded by $d_1^{-(N_1+1)}/2^{N_1+1}.$
Iteration will introduce another coefficient bounded by $C^{N_2}d_2^{-(N_2+1)}/2^{2(N_2+1)},$ then next by 
$C^{N_3}d_3^{-(N_3+1)}/2^{3(N_3+1)}.$ Since 
$$\frac{d_k^{-(N_k+1)}}{2^{k(N_k+1)}}=C^{N_k}\frac{(k^2)^{N_k+1}}{(2^k)^{N_k+1}},
$$
iteration at most $\log_2N$ times will lead to a product
$$\Pi_{k=1}^{\log_2N} C^{N_k}\left(\frac{k^2}{2^k}\right)^{\frac{N}{2^k}+1}\leq (C')^N$$
times a constant depending only on the first few derivatives of $u$ in the largest open set encountered. 

This yields the analyticity of $u$ in the smallest open set since all estimates are uniform in $N.$

\section{The case of non closed bicharacteristics with non trivial
  limit set}

\renewcommand{\theequation}{\thesection.\arabic{equation}}
\setcounter{equation}{0}
\setcounter{theorem}{0}
\setcounter{proposition}{0}
\setcounter{lemma}{0}
\setcounter{corollary}{0}
\setcounter{definition}{0}

We want to study a model of the form
$$ 
P (t, x, D_{t}, D_{x}) = D_{t}^{2} + X_2^{2}
$$
where
$$ 
X_{2} = g_{1}(x) g_{2}(x) [x_{1}D_{2} - x_{2} D_{1} + \mu (x_{1} D_{1} +
x_{2} D_{2}) + t^{k} (x_{1}D_{1} + x_{2}D_{2})],
$$
$$ 
g_{1}(x) = |x|^{2} - a^{2},
$$
$$ 
g_{2}(x) = b^{2} - |x|^{2},
$$
with $ 0 < a < b $ in the open set $ a < |x| < b $ and $ \mu > 0 $ is
a given constant.

The characteristic set of $ P $ in the above mentioned region is $
\Char (P) = \{\tau = 0, \ x_{1}\xi_{2} - x_{2}\xi_{1} + \mu (x_{1}
\xi_{1} + x_{2} \xi_{2}) + t^{k} (x_{1}
\xi_{1} + x_{2} \xi_{2})  = 0 \}$. 

As for the Poisson stratification of $ P $ we have
$$ 
\Sigma_{1} = \{\tau = 0, \ x_{1}\xi_{2} - x_{2}\xi_{1} + \mu (x_{1}
\xi_{1} + x_{2} \xi_{2}) + t^{k} (x_{1}
\xi_{1} + x_{2} \xi_{2})  = 0, \ t \neq 0 \};
$$
$$ 
\Sigma_{2} = \{\tau = t = 0, \ x_{1}\xi_{2} - x_{2}\xi_{1} + \mu (x_{1}
\xi_{1} + x_{2} \xi_{2}) = 0, x_{1} \xi_{1} + x_{2} \xi_{2} \neq 0 \};
$$
$$ 
\Sigma_{k+1} = \{0\},
$$
i.e. the zero section of $ R^{*}\R^{3} $ over the above specified
region. 

Evidently, since $ \codim \Sigma_{2} = 3 $, $ \Sigma_{2} $
(or rather its connected components) is not a symplectic submanifold
of $ R^{*}\R^{3} $.

Let us take a look at the Hamilton foliation of $ \Sigma_{2} $. 
Define
\begin{equation}
\label{eq:4.1}
A = \begin{bmatrix}
  \mu & 1 \\
  -1 & \mu
\end{bmatrix},
\end{equation}
so that
\begin{equation}
\label{eq:4.2}
\Sigma_{2} = \{ \tau = 0 = t, \ \langle x, A \xi \rangle = 0 \},
\end{equation}
where $ x = (x_{1}, x_{2}) $ and $ \xi = (\xi_{1}, \xi_{2}) $, $ \xi
\neq 0 $ and $ \langle x, \xi \rangle \neq 0 $. 

Then we know that $ \langle x, A \xi\rangle \equiv 0 $ on every leaf
in $ \Sigma_2 $, i.e. on every integral curve of the Hamilton field of 
$ \langle x, A \xi\rangle $ issued from a point in $ \Sigma_{2} $.

The Hamilton system is
\begin{eqnarray}
\label{eq:4.3}
\dot x & = & g_{1}(x) g_{2}(x) {}^{t}A x \\
\dot \xi & = & - g_{1}(x) g_{2}(x) A \xi. \nonumber
\end{eqnarray}
We easily see that, because of the structure of the matrix $ A $, we
have 
\begin{eqnarray*}
\frac{1}{2} d_{t} |x|^{2} & = & g_{1}(x) g_{2}(x) \frac{\mu}{2}
|x|^{2} \\
\frac{1}{2} d_{t} |\xi|^{2} & = & - g_{1}(x) g_{2}(x) \frac{\mu}{2} |\xi|^{2},
\end{eqnarray*}
so that both the spatial and the covariable  projections of the
bicharacteristics are logarithmic spirals. Moreover the spatial
projection spirals between the two asymptotic circles $ g_{i}(x) = 0
$, $ i = 1, 2 $, which are stationary orbits of the first two
equations in (\ref{eq:4.3}). 

We point out that $ d_{t} \langle x, \xi \rangle \equiv 0 $, so that $
\langle x, \xi \rangle $ is constant along the orbits and that once
the first two equations in (\ref{eq:4.3}) are solved the secon
couple---i. e. the covariable projection---is easy:
$$ 
\xi(t) = \exp \left[ - \int_0^{t} g_{1}(x(s)) g_{2}(x(s)) ds \ A \right]
\xi_{0}, 
$$
where $ \xi_{0} $ is its initial data.

We may apply to the operator $ P $ Theorem 4.2 in \cite{Sj-83} and
conclude that, if $ \gamma_{0} $ denotes a segment of a
bicharacteristic curve in $ \Sigma_{2} $, then either $ \gamma_{0}
\subset WF_{a}(u) $ or $ \gamma_{0} \cap WF_{a}(u) = \varnothing $,
where $ u $ is a solution of $ P u \in C^{\omega} $ in some open set.

\bigskip

Let now $ U $ be an open set in $ \R^{3} $ projecting onto an annulus
of the form $ a < |x| < b $ in the $ x $-variables. By iteratively
applying the above mentioned theorem one can prove the following

\begin{theorem}
\label{th:2}
Let $ u $ be a distribution such that $ P u \in C^{\omega}(U) $, $ U $
being defined as above. Then if both circles $ g_{i}(x) = 0 $ do not
intersect $ WF_{a}(u) $ we have that  $ u \in C^{\omega}(U) $.
\end{theorem}

\end{document}